\documentclass[12pt,reqno]{amsart}

\usepackage[colorlinks=false]{hyperref}
\usepackage{amssymb}
\usepackage{amsmath}
\usepackage{amsfonts}
\usepackage{float}
\usepackage{graphicx}
\usepackage{amsthm}
\usepackage{amscd}
\usepackage{calrsfs}
\usepackage{latexsym}
\usepackage{epsf}
\usepackage{oldlfont} 
\usepackage{xcolor} 

\graphicspath{{converted_graphics/}}

\begin{document}

\theoremstyle{plain}
\newtheorem{theorem}{Theorem}
\newtheorem{lemma}[theorem]{Lemma}

\theoremstyle{definition}
\newtheorem{definition}[theorem]{Definition}
\newtheorem{example}[theorem]{Example}
\newtheorem{conjecture}[theorem]{Conjecture}

\theoremstyle{remark}
\newtheorem{remark}[theorem]{Remark}

\newcommand{\seqnum}[1]{\href{https://oeis.org/#1}{#1}}  

\def \bN{{\mathbb N}}
\def \BB{{\cal B}}
\def \TT{{\cal T}}
\def \II{{\cal I}}

\title[A prime number ``game of life'']{A prime number ``game of life'': \\ can $\lfloor y\cdot p\#\rfloor$ be prime for all $p \geq 2$?}

\begin{abstract} \noindent
A new sequence in the spirit of the Mills primes is presented and its properties are investigated.
\end{abstract}

\author{Martin Raab}
\maketitle

\vspace{2mm}

\section{Introduction}

In the 1940s, W. H. Mills \cite{Mills} examined prime-representing functions, most notably he pointed out, based on a
result of Ingham \cite{Ingham2}, that there exists a constant $A$ such that $\lfloor A^{3^n} \rfloor$ is prime for all
integers $n \geq 1$. Mills' ansatz reproduces only a very sparse subset of primes. There are several variations to
obtain a larger subset of primes, e. g. by Elsholtz \cite{Elsholtz}, who conjectures the existence of a constant $A$
such that $\lfloor A^{(n+1)^2} \rfloor$ is prime for all $n \geq 1$; Kuipers \cite{Kuipers}, or Plouffe \cite{Plouffe},
who consider sequences of the form $\lfloor A^{c^n} \rfloor$ for $c$ as small as 1.01.

In this paper, we give an even more slowly growing formula---in fact, for its kind as an asymptotically
exponentially growing function, it exhibits the slowest possible growth rate while still producing prime numbers,
and a comparatively large subset at that.

While large portions of this paper consist of artisan work with well-known tools of number theory, there may
something to be gained in-between the lines, for example a way to calculate Chebyshev's $\theta(x)$ directly from
known values of $Li(x)-\pi(x)$ which we haven't seen before in literature. There are also a couple of open
problems which may or may not be worth investigating further.

\vspace{5mm}

\section{Definitions}
\vspace{2mm}

{\small

\noindent
\begin{tabular}{ll}
$p, q$ & variables for prime numbers \\
$p_n$ & the $n^{th}$ prime; ${p_n}\in \{2,3,5,7,11,13,17\ldots\}$ \\
$\pi(x)$ & the prime counting function; the number of primes $p_n\leq x$ \\
$p\#$ & the primorial function: $p\# = 2 \cdot 3 \cdot 5 \cdot 7 \cdot 11 \cdot \ldots \cdot p$ \\
& 1\# = 1 by convention \\
$\Phi(n,p)$ & number of integers $\leq n$ that are relatively prime to $p\#$ \\
$\log x$ & natural logarithm of $x$ \\
$\lfloor x\rfloor$ & the floor function; the largest integer $\leq x$ \\
$[a,b[$ & the open interval containing real numbers $x$ such that $a \leq x < b$ \\
$|x|$ & absolute value of $x$; $|x|=-x$ for $x<0$ \\
\end{tabular}
}

\newpage

\section{Setting the stage}

{\bf Lemma:}
Any real number $y$ such that $q = \lfloor y\cdot p\#\rfloor$ is a prime number for every $p \geq 2$
is contained in the interval $[1, 1.5[$. \\

If $y < 1$, then $\lfloor y\cdot 2\#\rfloor < 2$, which is less than the first prime number.

For $y \geq 1.5$, $\lfloor y\cdot 2\#\rfloor$ must be an odd prime $2m-1$ for some integer $m > 1$.
Then $\lfloor y\cdot 3\#\rfloor$ is a number in the interval $[6m-3, 6m-1]$, only the latter of which
can be prime because $6m-3$ is divisible by 3 and $6m-2$ is divisible by 2.
Likewise, any integer in the successive intervals $[p\# \cdot m-p, p\# \cdot m-1]$ is divisible by a
small prime $\leq p$ except for $p\# \cdot m-1$.

Hence there would have to be an integer $m$ for which $p\# \cdot m-1$ is always prime, i.e. the exponentially
growing sequence $f(0) = m-1$, $f(s+1) = (f(s)+1) \cdot p_s-1$ must only produce primes for all $s \geq 1$,
which is practically impossible (though we have not seen this being explicitly disproven yet!).

It is possible, however, to find arbitrarily long finite sequences for some $m$:
\begin{table}[H]
\caption{Smallest integers $m$ for which $p\# \cdot m-1$ is prime for $p \leq r$. The first unknown entry for
$r = 47$ has $m > 1.6 \cdot 10^{12}$ and is expected to occur when $10^{14} < m < 10^{15}$.
\label{table1}
}
\centering
\begin{tabular}{ccc|ccc}
\hline $r$ & $m$ & & & $r$ & $m$ \\
\hline
5 & 2 & & & 23 & 384427 \\
7 & 2 & & & 29 & 16114470 \\
11 & 9 & & & 31 & 259959472 \\
13 & 9 & & & 37 & 13543584514 \\
17 & 9 & & & 41 & 100318016379 \\
19 & 224719 & & & 43 & 100318016379 \\
\hline
\end{tabular}
\end{table}

\vspace{8mm}

To find out about an admissible value for $y$, we have to look for primes $q$, where each $q$ requires the
primality of $\lfloor q/p\# \cdot r\#\rfloor$ for every $r \leq p$. In the following setup, $Y$ will denote
the interval in which any number $y$ as described above can lie, setting $Y = [1, 1.5[$ for a start,
and then successively sharpen the bounds where, for every $p$, $Y \cdot p\#$ must contain at least one prime
to satisfy the conditions.

\newpage

\noindent
$\bullet\,$ Stage 1: $p = 2\,\rightarrow\, 2 \cdot Y = [2, 3[ \,\rightarrow\, \lfloor 2 \cdot Y\rfloor = 2$ -- no further action required yet \\
$\bullet\,$ Stage 2: $p = 3\,\rightarrow\, 6 \cdot Y = [6, 9[$ -- the only prime in here is 7, so Y narrows down to $[\frac{7}{6} , \frac{8}{6}[$ \\
$\bullet\,$ Stage 3: $p = 5\,\rightarrow\, 30 \cdot Y = [35, 40[$ -- one prime (37) in this interval, $Y$ shrinks to $[\frac{37}{30} , \frac{38}{30}[$ \\
$\bullet\,$ Stage 4: $p = 7\,\rightarrow\, 210 \cdot Y = [259, 266[$ -- once again, one prime in here and any possible $y$ lies between $\frac{263}{210}$ and $\frac{264}{210}$ \\
$\bullet\,$ Stage 5: $p = 11\,\rightarrow\, 2310 \cdot Y = [2893, 2904[$ -- now it’s getting interesting, since there are two primes in that interval, 2897 and 2903. We will split $Y$ into two separate intervals, $Y_1 = [\frac{2897}{2310} , \frac{2898}{2310}[$ and $Y_2 = [\frac{2903}{2310} , \frac{2904}{2310}[$ and proceed \\
$\bullet\,$ Stage 6: $p = 13\,\rightarrow\,$ each $Y_1$ and $Y_2$ brings forth one prime, 37663 and 37747 respectively. Cutting back both $Y$’s, then \\
$\bullet\,$ Stage 7: $p = 17\,\rightarrow\,$ $Y_1$ leads to 640279, whereas $Y_2$ gives 641701 and 641713, and $Y_2$ will be split into $Y_2$ and $Y_3$ in accordance with the primes \\
$\bullet\,$ Stage 8: $p = 19\,\rightarrow\, 9699690 \cdot Y_1$ contains a prime triplet, 12165311+$d$ for $d=\lbrace0,2,6\rbrace$, while $Y_2$ and $Y_3$ both fail to produce any further prime \\
$\bullet\,$ etc.

\vspace{8mm}

The results can be put into a nice grid:
\begin{table}[H]
\caption{This seems to be the new ``game of life'', a family tree of prime numbers.
Prime triplets $q+\lbrace0,2,6\rbrace$ or $q+\lbrace0,4,6\rbrace$ like in stage 8 are rather rare in the process,
one finds the next example at stage 152 ($p$ = 881), and there are no prime quadruplets $q+\lbrace0,2,6,8\rbrace$ in any stage $\leq$ 331.
\label{table2}
}
\centering
\begin{tabular}{c|c|c|c|c|c|c|c|c|c|c|c}
\hline
$\lfloor Y \rfloor$ & \multicolumn{11}{|c}{1} \\
\hline
$\cdot 2$ & \multicolumn{11}{|c}{+0} \\
\hline
$\cdot 3$ & \multicolumn{11}{|c}{+1} \\
\hline
$\cdot 5$ & \multicolumn{11}{|c}{+2} \\
\hline
$\cdot 7$ & \multicolumn{11}{|c}{+4} \\
\hline
$\cdot 11$ & \multicolumn{9}{|c|}{+4} & \multicolumn{2}{|c}{+10} \\
\hline
$\cdot 13$ & \multicolumn{9}{|c|}{+2} & \multicolumn{2}{|c}{+8} \\
\hline
$\cdot 17$ & \multicolumn{9}{|c|}{+8} & +2 & +14 \\
\hline
$\cdot 19$ & +10 & \multicolumn{6}{|c|}{+12} & \multicolumn{2}{|c|}{+16} & --- & --- \\
\hline
$\cdot 23$ & +16 & \multicolumn{6}{|c|}{+8} & +16 & +18 \\
\hline
$\cdot 29$ & +16 & \multicolumn{2}{|c|}{+6} & +14 & \multicolumn{2}{|c|}{+24} & +26 & +10 & --- \\
\hline
$\cdot 31$ & --- & +14 & +24 & --- & +6 & +14 & --- & --- \\
\hline
$\cdot ...$ & & +... & +... & & --- & +... \\
\hline
\end{tabular}
\end{table}

\newpage

Will the evolution go on forever? The following table shows how many partial intervals (or primes, respectively) $n$ there are after stage $s$, up to $s = 60$:
\begin{table}[H]
\caption{Talking about a bottleneck: the only prime at stage 29 is 350842542483891235293716663559065020274899073 ($\approx 3.5 \cdot 10^{44}$).
\label{table3}
}
\centering
\begin{tabular}{cccc|ccccc|cccc}
\hline $s$ & $p$ & $n$ & & & $s$ & $p$ & $n$ & & & $s$ & $p$ & $n$ \\
\hline
1 & 2 & 1 & & & 21 & 73 & 5 & & & 41 & 179 & 18 \\
2 & 3 & 1 & & & 22 & 79 & 4 & & & 42 & 181 & 17 \\
3 & 5 & 1 & & & 23 & 83 & 6 & & & 43 & 191 & 14 \\
4 & 7 & 1 & & & 24 & 89 & 3 & & & 44 & 193 & 24 \\
5 & 11 & 2 & & & 25 & 97 & 2 & & & 45 & 197 & 24 \\
6 & 13 & 2 & & & 26 & 101 & 1 & & & 46 & 199 & 28 \\
7 & 17 & 3 & & & 27 & 103 & 3 & & & 47 & 211 & 30 \\
8 & 19 & 3 & & & 28 & 107 & 1 & & & 48 & 223 & 36 \\
9 & 23 & 4 & & & 29 & 109 & 1 & & & 49 & 227 & 49 \\
10 & 29 & 6 & & & 30 & 113 & 3 & & & 50 & 229 & 44 \\
11 & 31 & 4 & & & 31 & 127 & 2 & & & 51 & 233 & 52 \\
12 & 37 & 5 & & & 32 & 131 & 5 & & & 52 & 239 & 53 \\
13 & 41 & 5 & & & 33 & 137 & 6 & & & 53 & 241 & 55 \\
14 & 43 & 9 & & & 34 & 139 & 12 & & & 54 & 251 & 67 \\
15 & 47 & 11 & & & 35 & 149 & 21 & & & 55 & 257 & 69 \\
16 & 53 & 10 & & & 36 & 151 & 19 & & & 56 & 263 & 72 \\
17 & 59 & 12 & & & 37 & 157 & 15 & & & 57 & 269 & 81 \\
18 & 61 & 8 & & & 38 & 163 & 16 & & & 58 & 271 & 79 \\
19 & 67 & 6 & & & 39 & 167 & 24 & & & 59 & 277 & 85 \\
20 & 71 & 11 & & & 40 & 173 & 18 & & & 60 & 281 & 83 \\
\hline
\end{tabular}
\end{table}

Although the population increases considerably after stage 33, the data doesn’t provide too much confidence
that it will continue to do so. We have to take a closer look at why this sequence is so erratic.

\section{Analyzing the game}

The size of the primes at each stage $s$ with the corresponding $p$, $p_s$ being the $s^{th}$ prime,
is always about $1.25 \cdot p_s\#$. The probability of a random number of this size being prime
is $1/(\log p_s\# + 0.23)$, or $1/(\theta(p_s) + 0.23)$ where $\theta(p)$ is the first Chebyshev function.
For one prime $q$ of the stage $s-1$ there has to be on average at least one prime in the interval
$[q \cdot p_s, (q+1) \cdot p_s]$ to have a chance that the sequence keeps on producing ever more primes.
More precisely, the probability of getting at least one prime out of the respective interval
must be bigger than getting no prime at all.
\begin{equation} \label{eq1}
\frac{\binom{p-1}1}{(\log q - 1)^1}(1-\frac{1}{\log q})^{p-1} > \frac{\binom{p-1}0}{(\log q - 1)^0}(1-\frac{1}{\log q})^{p-1} \: \text{ or simply }
\: \frac {p-1}{\theta(p)-0.77}>1.
\end{equation}\
This is true for $p \geq 3$ if and only if $\theta(p) < p-0.23$, or $\theta(p) < p$ for short,
since the constant value $0.23 = \log y$ is negligible for large $p$. \\

Now $\theta(p)$ is in fact smaller than $p$ most of the time --- but only slightly so. And indeed $p-\theta(p)$
behaves just like the mercurial $Li(p)-\pi(p)$, along with the infinitude of sign changes (see Littlewood \cite{Littlewood}): \\

Using the simple sum \footnote{
~As opposed to the ``European convention'' $Li(x) = \int_2^p \frac{\text{d} x}{\log x}$ and the ``American convention''
$Li(x) = \int_0^p \frac{\text{d} x}{\log x}$ which each differ from said sum by less than 1 for $p \geq 12.0050107...$ (European) or $p \geq 2$ (American).
}
\begin{equation} \label{eq2}
Li^*(p) = \sum_{x=2}^p \frac{1}{\log x} \: ,
\end{equation}\
then
\begin{equation} \label{eq3}
p-\theta(p) = (\log p) (Li^*(p)-\pi(p))+1- \sum_{x=2}^{p-1} \bigg(\log (\frac{x+1}{x}) (Li^*(x)-\pi(x))\bigg).
\end{equation}\

This formula makes it also clear that $p-\theta(p)$ changes sign quite a while before $Li(p)-\pi(p)$ acquires negative values.
For the sake of a ballpark figure: assuming $Li(x)-\pi(x)$ is on average close enough to $\frac{\sqrt{x}}{\log x}$,
then $p-\theta(p)$ drops below zero by the time $Li(p)-\pi(p) < \frac{2 \sqrt{p}+O(p^{1/3})} {(\log p)(\log p -2)}$.
It can be expected that this happens for the first time close to the first $Li(x)-\pi(x)$ crossover near $1.4\cdot 10^{316}$
(Bays and Hudson \cite{Bays-Hudson}, more extensive calculations by Saouter and Demichel \cite{Saouter-Demichel}).
The results of Platt and Trudgian \cite{Platt-Trudgian} also confirm this. \\
\\

Empirical data and heuristical reasoning suggests that $p-\theta(p)$ can usually be found in the vicinity of $\sqrt{p}$.
The bias is given by Riemann’s prime counting formula \cite{Riemann}
\begin{equation} \label{eq4}
\pi(x) = \sum_{n=1}^\infty \frac{\mu(n)}{n} \bigg(Li(x^{1/n})-\sum_\rho Li(x^{\rho/n})+\int_{x^{1/n}}^\infty \frac{\text{d}u}{u(u^2-1)\log u}-\log(2)\bigg).
\end{equation}\

If the Riemann Hypothesis (RH) is true, $Li(p)-\pi(p)$ oscillates along $\frac{\sqrt{p}}{\log p}$ with a magnitude of at most
$O(\sqrt{p} \log p)$ (von Koch \cite{Koch}), and so $p-\theta(p)$ oscillates along $\sqrt{p}$ with a magnitude of at most
$O(\sqrt{p} \log^2 p)$. Ingham \cite{Ingham} also proved this bias for sufficiently large $p$. \\

On a large scale, the number of intervals should increase slowly and asymptotically with $p-\theta(p)$,
thus it may also decrease whenever $\theta(p) > p$. Yet the actual development strongly depends
on the local irregularities of the distribution of prime numbers, so one can only assign the probability
that the initial conditions are adequate. \\

Admittedly, \eqref{eq1} wasn’t quite accurate. The number of primes in an interval immediately after a semiprime as it is
the case, especially in this short type of interval, is on average a bit smaller than for a randomly chosen interval.
The following table illustrates this sort of difference for a small interval of ten numbers:
\begin{table}[H]
\caption{For a random interval start and interval length 10, the average of numbers that are coprime to $2310 = 11\#$
is above two, while for an interval of the same length after a prime number (or any number coprime to $11\#$), significantly
fewer than two such numbers are to be expected. Hardy and Littlewood \cite{Hardy-Littlewood} laid a solid foundation regarding this issue.
\label{table4}
}
\centering
\begin{tabular}{c|c|ccccc|c|cc}
\hline
& & \multicolumn{5}{|c|}{divisible by} & & not & interval \\
number & & 2 & 3 & 5 & 7 & 11 & & divisible & sum \\
\hline
random +[1..10] & & 50\% & 33\% & 20\% & 14\% & 9\% & & 21\% & \bf{2.08} \\
\hline
prime +1 & & 100\% & 50\% & 25\% & 17\% & 10\% & & 0\% & 0 \\
prime +2 & & 0\% & 50\% & 25\% & 17\% & 10\% & & 28\% & 0.28 \\
prime +3 & & 100\% & 0\% & 25\% & 17\% & 10\% & & 0\% & 0.28 \\
prime +4 & & 0\% & 50\% & 25\% & 17\% & 10\% & & 28\% & 0.56 \\
prime +5 & & 100\% & 50\% & 0\% & 17\% & 10\% & & 0\% & 0.56 \\
prime +6 & & 0\% & 0\% & 25\% & 17\% & 10\% & & 56\% & 1.13 \\
prime +7 & & 100\% & 50\% & 25\% & 0\% & 10\% & & 0\% & 1.13 \\
prime +8 & & 0\% & 50\% & 25\% & 17\% & 10\% & & 28\% & 1.41 \\
prime +9 & & 100\% & 0\% & 25\% & 17\% & 10\% & & 0\% & 1.41 \\
prime +10 & & 0\% & 50\% & 0\% & 17\% & 10\% & & 38\% & \bf{1.78} \\
\hline
\end{tabular}
\end{table}

Ultimately, the number of numbers that are relatively prime to $k\#$ in the intervals in question $[q\cdot p+1, (q+1)\cdot p-1]$
--- denoted by $\Phi(p-1, k)$ --- is
\begin{eqnarray} \label{eq5}
\begin{split}
\prod\limits_{\substack{u=3\\ u\,prime}}^k (1-\frac{1}{u-1})\cdot\bigg(\lfloor\frac{p-1}{2}\rfloor
+\sum\limits_{\substack{v=3\\v\,prime}}^{\min(\frac{p-1}{2},k)}\frac{\lfloor\frac{p-1}{2v}\rfloor}{v-2}
+\sum\limits_{\substack{v_1=3\\v_1\,prime}}^{\min(\sqrt{\frac{p-1}{2}},k)} \sum\limits_{\substack{v_2=v_1+2\\v_2\,prime}}^{\min(\frac{p-1}{2v_1},k)}\frac{\lfloor\frac{p-1}{2v_1v_2}\rfloor}{(v_1-2)(v_2-2)} \\
+\sum\limits_{\substack{v_1=3\\v_1\,prime}}^{\min(\sqrt[3]{\frac{p-1}{2}},k)}
 \sum\limits_{\substack{v_2=v_1+2\\v_2\,prime}}^{\min(\sqrt{\frac{p-1}{2v_1}},k)} \sum\limits_{\substack{v_3=v_2+2\\v_3\,prime}}^{\min(\frac{p-1}{2v_1v_2},k)}\frac{\lfloor\frac{p-1}{2v_1v_2}\rfloor}{(v_1-2)(v_2-2)(v_3-2)}+...\bigg).
\end{split}
\end{eqnarray}\

Yet disregarding that $p$ itself doesn’t appear as a factor in said interval, this particular formula is only valid for $2<k<p$.
In contrast, for a random interval $R$, $\Phi_R(p-1, k) = (p-1)\cdot W(k)+O(1)$, where
\begin{equation} \label{eq6}
W(k) = \prod\limits_{\substack{u=2\\u\,prime}}^k \bigg(1-\frac1u\bigg).
\end{equation}\

If $\Phi(p-1, k)$ is then divided by $W(k)\cdot\frac{p-2}{p-1}$ (as $p$ itself doesn't appear as a factor in the intervals
we're looking at), then the result is an ``adjusted'' interval length.
For this, $k$ must be appropriately large to attain the desired value ($k>\log p$, say).

There is a special connection between said adjusted interval length and the twin prime constant $C_2=0.6601618…$:
for $k\to\infty$, the adjusted interval length $\Phi(p-1, k)/W(k)$, on the basis of a prime number preceding the
interval, can be expressed as $C_2\cdot(p-1+a)$, with $a\in\mathbb Q$ being defined by
\begin{equation} \label{eq7}
a=2\cdot \sum_{x=1}^\frac{p-1}{2} \bigg(\prod\limits_{\substack{r:\\\text{every distinct}\\\text{odd prime}\\\text{factor of }x}} \frac{r-1}{r-2}-1\bigg).
\end{equation}\

Some values of $a$ include
\begin{table}[H]
\caption{The value of $\frac{a}{p-1}$ is asymptotic to $\frac1{C_2}-1-O(\frac{\log p}{2p\cdot C_2})$.
\label{table5}
}
\centering
\begin{tabular}{ccc|ccc}
\hline $p-1$ & $a$ & & & $p-1$ & $a$ \\
\hline
2 & 0 & & & 20 & 116/15 \\
4 & 0 & & & 30 & 6866/495 \\
6 & 2 & & & 40 & 141274/8415 \\
8 & 2 & & & 50 & 1329632/58905 \\
10 & 8/3 & & & 100 &  129132288244/2731483755 \\
12 & 14/3 & & & 150 & 123443421975532168/1666745013838905 \\
\hline
\end{tabular}
\end{table}

As $k\to\infty$, the adjusted interval length (for our case with the semiprimes) is --- apart from $p = 3, 7,\text{ and }13$
--- always a bit smaller than $p-1$. Dividing it again by $\log(q\cdot p)$ and we arrive at the expected average number of
primes in one interval. The ratio of the adjusted interval length vs. $p-1$, which is equivalent to
$\frac{\Phi(p-1, k)}{W(k)\cdot(p-2)}$ for appropriately large $k$, will be denoted below by $\psi(p)$. \\

Here we were focusing on heuristical aspects. For further reading about the structural analysis of the evolutional
nature of the sequence, and potential vistas, we refer the reader to Santana \cite{Santana} who gives an introduction
on the notion of evolution on sets.

\section{Proposing practical predictions}

Using these heuristics, which are so far in very good agreement with the factual data, we can start to calculate
the probabilities for the game to continue. \\

For example, there are 594 primes at stage 100, where $p = 541$ and every prime $q=\lfloor p\#\cdot Y\rfloor$
has 220 decimal digits. Each $q$ then has a certain chance to yield $n$ primes in the following stages:

\begin{table}[H]
\caption{Only the fittest will survive.
\label{table6}
}
\centering
\begin{tabular}{ccccccccccc}
\hline $s$ & $p$ & $n=0$ & $n=1$ & $n=2$ & $n=3$ & $n=4$ & $n=5$ & $n=6$ & $n=7$ & $n>7$  \\
\hline
101 & 547 & 34.52\% & 36.75\% & 19.53\% & 6.90\% & 1.83\% & 0.39\% & 0.07\% & 0.01\% & 0.00\% \\
102 & 557 & 49.78\% & 19.46\% & 14.17\% & 8.22\% & 4.33\% & 2.16\% & 1.03\% & 0.47\% & 0.36\% \\
103 & 563 & 58.57\% & 12.14\% & 10.10\% & 7.05\% & 4.63\% & 2.95\% & 1.83\% & 1.12\% & 1.62\% \\
104 & 569 & 64.34\% & 8.33\% & 7.46\% & 5.75\% & 4.22\% & 3.03\% & 2.14\% & 1.50\% & 3.22\% \\
105 & 571 & 68.45\% & 6.11\% & 5.72\% & 4.68\% & 3.68\% & 2.84\% & 2.17\% & 1.64\% & 4.71\% \\
... & ... & ... & ... & ... & ... & ... & ... & ... & ... & ... \\
110 & 601 & 78.69\% & 2.06\% & 2.14\% & 1.99\% & 1.80\% & 1.62\% & 1.44\% & 1.27\% & 8.99\% \\
150 & 863 & 90.06\% & 0.10\% & 0.11\% & 0.11\% & 0.11\% & 0.11\% & 0.11\% & 0.11\% & 9.20\% \\
200 & 1223 & 91.20\% & 0.01\% & 0.02\% & 0.02\% & 0.02\% & 0.02\% & 0.02\% & 0.02\% & 8.69\% \\
1000 & 7919 & 91.48\% & 0.00\% & 0.00\% & 0.00\% & 0.00\% & 0.00\% & 0.00\% & 0.00\% & 8.52\% \\
2000 & 17389 & 91.48\% & 0.00\% & 0.00\% & 0.00\% & 0.00\% & 0.00\% & 0.00\% & 0.00\% & 8.52\% \\
\hline
\end{tabular}
\end{table}

The column for $n = 0$ is of importance to find out about the overall probability. (Once an interval is devoid
of primes, it never comes back!) Within the first few thousand values of $s$,
this value approaches $\omega\approx 91.476\%$. This gives the overall probability of $\omega^{594}\approx
1.039\cdot 10^{-23}$ that the sequence fails at this point.

Generally the value $\omega(s)$ can be defined for every $s\geq 1$ as the probability that a prime $q$ at stage $s$
will not cause an infinite sequence of primes in the further process.

For convenience, a recursive formula can be considered where, with $\psi(p) = \frac{\Phi(p-1, k)}{W(k)\cdot(p-2)}$ from above,

\begin{equation} \label{eq8}
\omega(s-1) = \bigg(1-\frac{\psi(p)}{\log q}\bigg)^{p-1} \sum_{j=0}^{p-1} \frac{\binom{p-1}j \omega(s)^j}{(\frac{\log q}{\psi(p)}-1)^j}
\end{equation}\\
and a starting point $\omega(t) = 1-1/\sqrt{2p_t}$ for some $t$ can be set to quickly compute the values for all $s < t$.
Varying $\omega(t)$ between 0 and $1-\varepsilon$ for some $\varepsilon > 0$ leads to a lower and upper bound on $\omega(s)$.
It may be useful that if $\omega(t) = \omega(s)+\delta$ then $\omega(t+1) \sim \omega(s+1)+\delta\cdot(1-\omega(s+1))/2$. \\

Again, we should bear in mind that these primes are not as flexible as \eqref{eq8} takes them to be:
the probability of obtaining $n$ primes in the following step is zero when $n$ exceeds the maximal possible number
of primes in the given interval (e.g. $n > 97$ for $p-1 = 546$) as fixed in the specific rules for the k-tuple conjecture.
To take account of that --- at least to some degree ---, $W(k)$ from \eqref{eq6} is deployed again such that $k\# < p$.
For $p = 547$, we set $k = 7$ and $W(k)$ = 48/210 = 1/4.375. As $p$ gets larger than 11\# = 2310, $W(k)$ can be
changed to 480/2310 = 1/4.8125 and so on.

\begin{equation} \label{eq9}
\omega(s-1,k) = \bigg(1-\frac{\psi(p)}{W(k) \log q}\bigg)^{W(k)\cdot(p-1)} \sum_{j=0}^{W(k)\cdot(p-1)} \frac{\binom{W(k)\cdot(p-1)}j \omega(s,k)^j}{(\frac{W(k) \log q}{\psi(p)}-1)^j}
\end{equation}\

This refined calculation reveals a slightly smaller value for $\omega(100,7)$ (compared to $\omega(100)$), namely 91.435\%.
The entire probability that the sequence fails at this point drops by 24\% to $7.937\cdot 10^{-24}$. To show how vague these
percentages are (especially at this rather early stage of computation), the actual portion of primes that don’t survive after stage 100 ---
at least 547 out of 594 --- is $\geq 92.088\%$; to the $594^{th}$ power, this is over 68 times more than the value from the
``refined'' calculation! In other terms, 47 surviving primes as opposed to a predicted number of 51, from this point
of view the deviation is not that bad (and it gets better for larger $p$). \\
\\
$\omega(100) = \omega(100,1) = 91.4760031\%$ \\
$\omega(100,2) = 91.4637729\%$ \\
$\omega(100,3) = 91.4515074\%$ \\
$\omega(100,5) = 91.4422850\%$ \\
$\omega(100,7) = 91.4345844\%$ \\

It is difficult to pin down exact values there, but the refined calculation should at least give some indication
on how large the error might be. \\

Not too surprisingly, $\omega(s)$ is largely situated in the neighborhood of $1-1/\sqrt{2p_s}$. This can be verified heuristically as follows:
Proceed as to be seen in table \ref{table6}, starting with stage $s$. Consider regular intervals with evenly distributed
primes $q$, the size of the intervals being $z\cdot\log(q)$ for a chosen $z\geq 1$, so for lim $q\to\infty$, the probability
that one interval doesn’t contain a prime at stage $s+1$ is $e^{-z}$.
The probability to get 0 primes after $t$ steps at stage $s+t$ is $f(s+t)$, where $f(s) = 0$ and recursively $f(x+1) = e^{z\cdot(f(x)-1)}$.
As $\lim z\to 1$ for $s$ fixed, $x = 1-2(z-1)+O((z-1)^2)$, which can be approximated by $1-1/\sqrt{2p_s}$ for
$z = p/\theta(p) \approx 1+1/\sqrt{p}$ (on average). \\

What is more, for $\theta(p) = p+m\cdot \sqrt{p}$, where $m = -1+O(log^2p)$ (under the RH),

\begin{equation} \label{eq10}
\omega(s+1)-\omega(s) \sim (\omega(s)-1) \bigg(\frac{m}{\sqrt{p}}+\frac1p+\frac{1-\omega(s)}2+O\bigg(\frac{\log p}{\sqrt{p^3}}\bigg)+O\bigg(\frac{1-\omega(s)}{\sqrt{p}}\bigg)\bigg)
\end{equation}\

with error terms depending on either $p$ or $\omega(s)$. As long as $m\approx -1$ as is expected on average in the long run,
if $\omega(s)$ is chosen a bit larger than $1-1/\sqrt{2p_s}$, the values for the subsequent $\omega(s+x)$ by the recursion formula
above quickly head toward 1, so the first term in the parenthesis above, $m/\sqrt{p}$, becomes the most significant.
If then at some point $Li(p)\approx\pi(p)$, thus $m\approx 0$, while $\omega(s)$ is very close to one, then the second term 1/$p$
becomes the most significant, but doesn’t nearly have the (then opposite) impact on $\omega(s)$ as in the region where $m\approx -1$.
Now what does that say? If $\omega(s)$ stays close to 1, the initially chosen value is larger than the actual value, which is good because
the actual value --- the probability that the sequence terminates --- is desired to be as small as possible. An $m$-value of +1
would adequately counter the effect of the expected average value $-1$, then again that should happen just as often as $m \leq -3$.
It should be noted that there is yet no known effective region where $|m+1|\geq 2$, making it a candidate for a future endeavor.
\footnote{~Myerscough \cite{Myerscough} tabulates the logarithmic density for the case $m > 1$ as $1.603 \cdot 10^{-95}$.} \\

Once the number of primes $n$ --- the number of possibilities for $y$ --- reaches a ``stable'' value, one may give a rough estimate
as to how the sequence continues. For instance, with $p$ = 1153 we have $n$ = 19690. In the next stage, with $p$ = 1163, one is
likely to find approximately $n = 19690\cdot\psi(p)\cdot(p-1)/(\theta(p)+0.23)\approx 20298$ primes. \\
Most of the time, the actual value lies between $n-\sqrt{n}$ and $n+\sqrt{n}$, and, by standard deviation measures, more than 99.7\%
of the time it should be well within $n\pm 3\sqrt{n}$. But even considering a negative deviation of $4\sqrt{n}$: $-4\cdot\sqrt{20298} = -569.88...$
is still outnumbered by the predicted growth $20298-19690 = 608$, corroborating that the initial assumption holds with a very high
probability if the distribution is normal (a sufficient but not particularly necessary condition).

\begin{figure}[H]
  \centering
  \includegraphics[bb=-40 0 900 600,width=7.5in,height=5in,keepaspectratio]{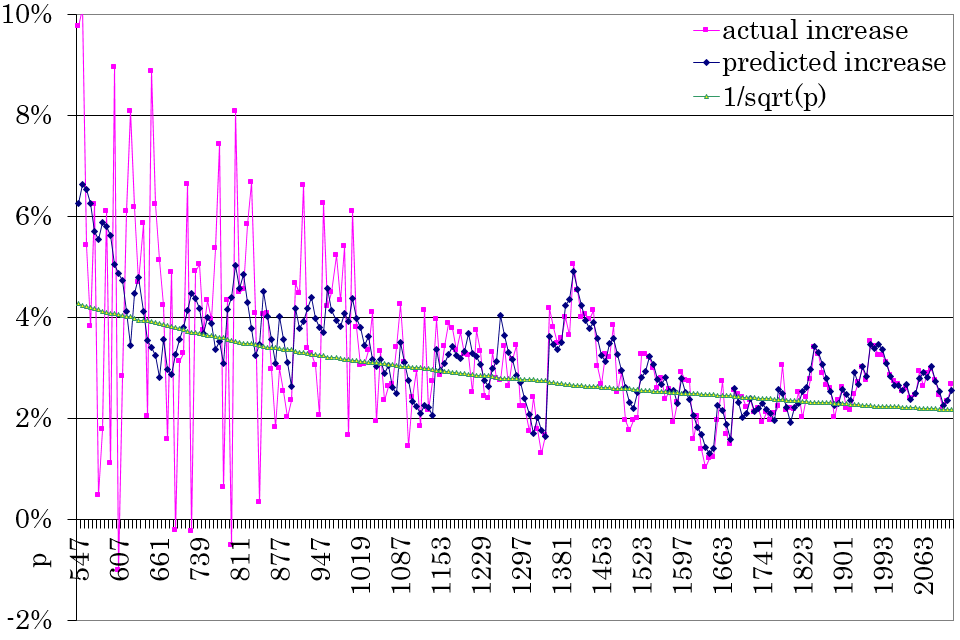}
  \caption{Comparison. With the prime number output becoming more generous, the forecast is increasingly accurate.}
\label{fig1}
\end{figure}

Expanding this supposition, we can say the sequence is ``stable'' when $n > p$, where the average absolute growth rate
$n/\sqrt{p}$ is bigger than the majority of the local deviations of order $\sqrt{n}$. This critical boundary is
firmly overstepped at $s = 98$: $p = 521$ and $n = 559$ (compared to $s = 97$ with $p = 509$ and $n = 490$). \\

Looking at predicted vs. actual increase in terms of standard deviations, so far the sequence behaves pretty ``normal'':

\begin{figure}[H]
  \centering
  \includegraphics[bb=-100 0 900 600,width=6in,height=4in,keepaspectratio]{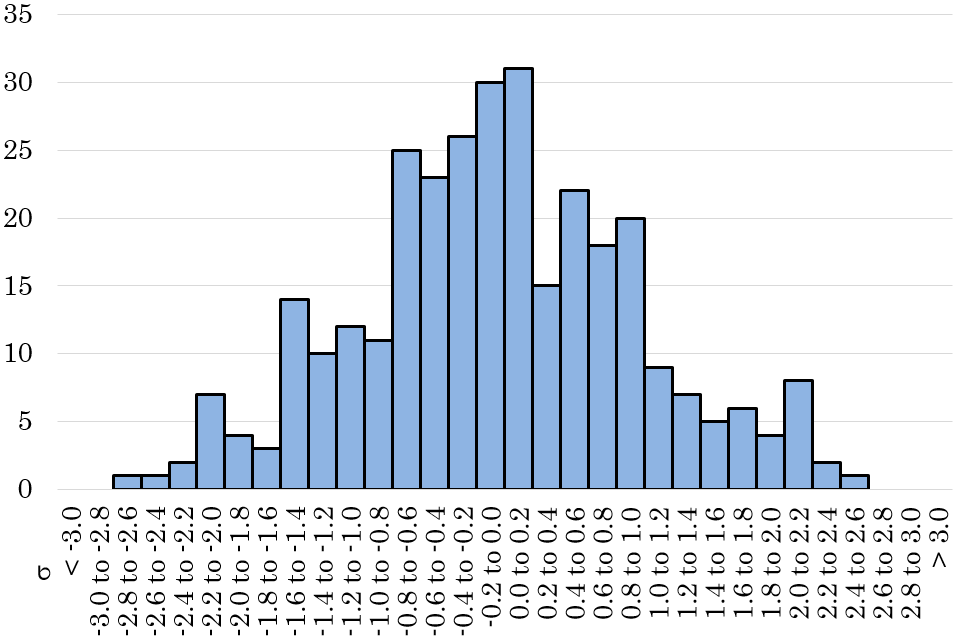}
  \vspace{2mm}
  \caption{Histogram of the standard deviations from the predicted number of primes at each stage $s = 99…251$.
70.6\% of the values are within 1$\sigma$, 92.2\% within 2$\sigma$, and all within 3$\sigma$.
There are 78 negative deviations vs. 75 positive ones.}
\label{fig2}
\end{figure}

Long-term predictions imply that the number of primes may surpass $10^6$ at $\approx$ stage 339 ($p = 2281$),
$10^9$ at $\approx$ stage 702 ($p = 5297$), $10^{12}$ at $\approx$ stage 1210 ($p = 9811$), $10^{100}$ at
$\approx$ stage 118948 ($p = 1568341$), and may have reached $n \approx 10^{10^{155}}$ by the time that
$\theta(p) > p$ where we expect to see the next slight decline of the sequence ever since stage 139 (see table \ref{table7}). \\

The following seemingly unrelated approximation formula for $n$ can be used in the medium term:

\begin{equation} \label{eq11}
n \approx \frac{\sqrt{e^{\sqrt{3s}}}}8
\end{equation}\\
which is surprisingly accurate in spite of its simplicity in the narrow sense that it is off by less than 20, with
one little exception, for $s \leq 61$, off by at most 20\% for $98 \leq s \leq 114$, and off by at most 15.4\% for
$115 \leq s \leq 331$ --- compared to more elaborate predictions, the latter may hold even for a couple more $s \leq 500$.
Furthermore, \eqref{eq11} may be off by a factor of no more than 2 for $66 \leq s \leq 1000+$. In the long run, however,
the expression in \eqref{eq11} is expected to far exceed $n$, thus calling into question any more daring predictions based on it.

\vspace{5mm}

\section{The status quo, pt. I: risk of failure} 

At stage 331 ($p = 2221$) there are 810387 probable primes, and $\omega(331) \approx 95.45431\%$.
$0.9545431^{810387} \approx 3.71 \cdot 10^{-16374}$ --- that is, for the time being, the conjectured probability that
the game fails for some $s > 331$.

Up to that point, 32103994 numbers have not yet been tested for certified primality, the smallest of those being in the
ballpark of $2.38 \cdot 10^{247}$. No counterexample is known for the combination of a Baillie-Pomerance-Selfridge-Wagstaff
pseudoprimality test and a strong Lucas test which were carried out on the probable primes.
\footnote{~including some trial division, additionally decreasing the chances that a pseudoprime slipped through}
However if there is one among all the established PRPs, it could take effect on at most 85448 PRPs (the most prolific branch at $p = 601$).
$0.9545431^{810387-85448} \approx 1.00 \cdot 10^{-14647}$, i.e. $2.70 \cdot 10^{1726}$ times larger. It is yet unclear
how much this will affect the above evaluation, especially by taking a look at Pomerance \cite{Pomerance} who gives a
heuristic argument that the number of counterexamples for a BPSW primality test up to $x$ is $\gg x^{1-\varepsilon}$ for
some $\varepsilon > 0$ and sufficiently large $x$. In either case, said ratio can be systematically reduced by checking
the primality of the prolific branches in the sequence. 

\vspace{5mm}

\section{The status quo, pt. II: closing in on $y$}

Now, what about the initially proposed values for $y$? The lower and upper bounds are easily obtained by picking the
smallest prime $q_{s,1}$ and the largest prime $q_{s,n}$ of the last computed stage and dividing by $p\#$:
$y_{\min} = q_{s,1}/p\#$ and $y_{\max} = (q_{s,n}+1)/p\#$.
The calculations so far demonstrate that any $y$ lies in a very small range between
$$
\mbox{\small 1.2541961015780119362776795549142134237798692180426221958327225546088646994287514475}
$$
$$
\mbox{and}
$$
$$
\mbox{\small 1.2541961015780119362776795549142134237798692180426221958327225546088646994290445894.}
$$
In particular, $y_{\max} - y_{\min} = 2.931419... \cdot 10^{-76}$. \\

Although it is not dead certain that such a $y$ exists, if it actually does, then there should be infinitely many of those.
Otherwise, we’d have to assume a unique solution where the probability of getting one prime for $p \to \infty$ --- the column
for $n = 1$ in Table \ref{table6} --- is greater than zero. But if this is the case, then at every stage a factor of
approximately $1/e$ falls back on the column for $n = 0$ which then would converge to 1, challenging our assumption that
there is a unique solution. It would be nice, however, to have a rigorous proof. \\

\section{Ancestors and survivors}

Having calculated all primes up to a certain stage, we can trace back these primes for previous stages and see which of those
have ``survived'' in the long run. All probable primes from stage 331 originated from one prime at stage 43.
Recall that there were a total of 14 primes at stage 43, so the other 13 ``petered out'' along the way, including one branch
that originated in stage 37 tenaciously keeping up until it succumbed at stage 95.

The first split then appears to be at stage 44. From there on, at least two values of $y$ may satisfy the hitherto
harsh conditions. And suddenly there’s a lot more to come:

\begin{figure}[H]
  \centering
  \includegraphics[bb=-160 0 1008 672,width=7in,height=5in,keepaspectratio]{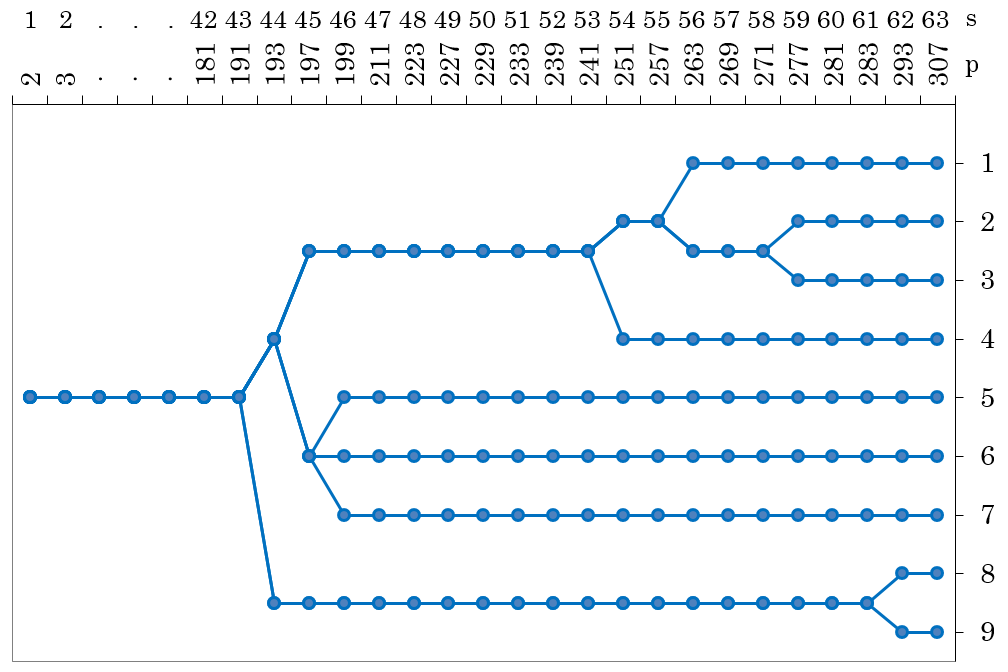}
  \caption{Number of possible values for $y$ (for all we know). From $s = 95$ onward, all primes 
$q$ emerge from one common ancestor at $s = 43$, namely 
{\small 1292942159746921794791923187781692727375711831825607985864285936838920812561.}
Naturally, the branches in the picture are not equally strong. For $s = 331$ ($p = 2221$), the 
following number of probable primes divide up on the nine given groups:
11477 + 9305 + 64706 + 123386 + 137475 + 101078 + 13648 + 77384 + 271928.}
\label{fig3}
\end{figure}

If one was to specify these values as $y_0$, $y_1$ and so on to attain fixed constants for the various 
possibilities, one might set $y_0$ as the minimum value $y_{min}$ as described above, and label each 
of the consecutive numbers $y_1$, $y_2$ … as the other smallest possible value(s) after a 
precalculated split, when it is ``safe enough to assume'' that those branches are stable.

So $y_0 = y_{min}$ (see above), and then \\
$y_1 = y_0 + 2.9314187260027917698243… \cdot 10^{-76}$ \\
(from the first split at stage 44, $\lfloor y_0 \cdot 191\#\rfloor \cdot 193+30$ and $+88$, $y_1 \approx y_0+58/193\#$) \\
$y_2 = y_0 + 1.9278096567063412938136… \cdot 10^{-79}$ \\
$y_3 = y_2 + 1.1844960407807006213569… \cdot 10^{-80}$ \\
$y_4 = y_2 + 1.5406827193677358123185… \cdot 10^{-80}$ \\
$y_5 = y_0 + 1.5619786064199049976025… \cdot 10^{-100}$ \\
$y_6 = y_0 + 1.4764435119323593401152… \cdot 10^{-104}$ \\
$y_7 = y_6 + 1.4210514603041973809966… \cdot 10^{-111}$ \\
$y_8 = y_1 + 2.5503871104036223556538… \cdot 10^{-119}$ \\
… \\

The value $y_6$ is one that may be doubted more than $y_0$ ... $y_5$ since it derives from the weakest 
branch of Figure 3 with only 9305 primes. But since 9305 $>$ 2221 (thus a stable branch), it can be expected to 
remain valid. \\

The stability criterion as described above, $\lim_{sup} n/p < 1$ for finite branches, may even be strengthened to
$\lim_{sup} n/({\frac12} \cdot \sqrt p \cdot \log p) < 1+\varepsilon$ and still hold for the entirety of the sequence for
relatively small $\varepsilon$. Even if this is not true for all $\varepsilon > 0$ (and under certain conditions ---
looking at runs of consecutive non-surviving primes of a given stage --- can be anticipated to fail),
in practice, it will be a safe bet that no-one will explicitly find a counterexample for
branches with more than $\frac12 \cdot \sqrt p \cdot \log p$ primes.
There might be a bound $\lim_{sup} n/({\frac12} \cdot \sqrt p \cdot \log p) \leq c$ for finite branches for some constant
$c$, and thus far $c > 0.75478$ ($177^{th}$ branch/prime of stage 111 with 64 primes at stage 123 but
dissipating at stage 316). So two more questions are probably to remain unanswered: Is
there such a bound, maybe with $c > 1$? If yes, what is the heuristically/actually largest $c$? If
no, is there any other bound?

\vspace{8mm}

The bifurcation tree as depicted in Figure \ref{fig3} gets quite impressive when expanded. For example, tracing back the
primes from stage 331 to stage 100 crystallizes out 47 primes (compared to originally 594) with an eventful history.

\begin{figure}[H]
  \centering
  \includegraphics[bb=10 0 1008 672,width=9in,height=6in,keepaspectratio]{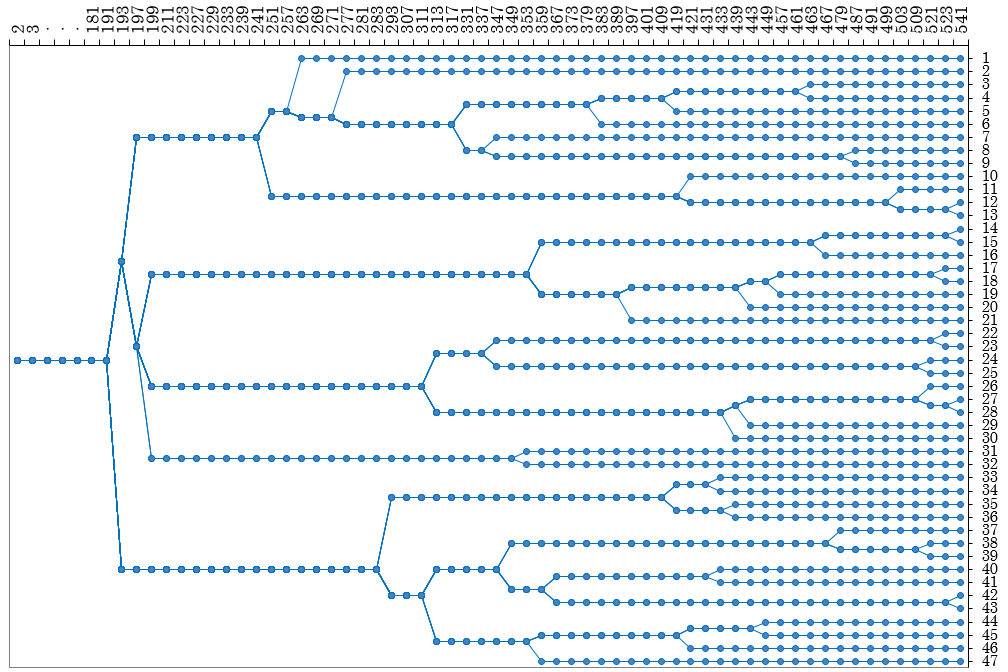}
  \caption{Alleged number of possibilities up to stage 100.}
\label{fig4}
\end{figure}

It is interesting to note the threefold split at stage 46, $p = 199$. If all of these primes persist an infinite run,
this still is the only point for $s \leq 100$ where more than two primes out of a single prime from the previous stage
will simultaneously continue the race. But how probable is this now? We have to consider the possible number of
primes arising from one prime of the previous stage and the possibility that at least three of them become stable
sequences, including the permutations if more than three primes show up simultaneously in the current stage:

\begin{equation} \label{eq12}
\sum_{x=3}^\infty \frac{\binom{p-1}x}{(\log q -1)^x} (1-\frac 1{\log q})^{p-1} (1-\omega(s-1))^x \binom x3
\end{equation}\

For $s = 101$, this value is 0.0043\%, thus $\approx$ 2.76\% for all 594 primes of stage 100.

\begin{figure}[H]
  \centering
  \includegraphics[bb=-40 0 840 560,width=7.5in,height=5in,keepaspectratio]{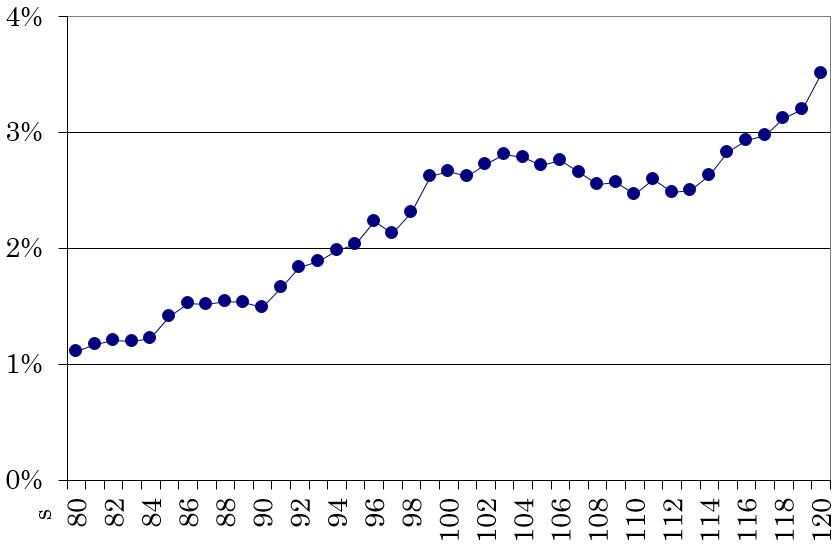}
  \caption{The value is growing in the process (e.g. about 52\% for $s = 252$), so we can soon expect another threefold split.
In fact, the next one appears to be in the upper region of $p = 613$, $s = 112$.}
\label{fig5}
\end{figure}

Substituting 3 for $k$ in \eqref{eq12}, it can handily be used for the probability of a $k$-fold split for any $k \geq 2$.

The possibly last stage without any split, that is without a temporarily increasing number of surviving primes and
hence admissible values for $y$, is found at $s = 103$ ($p = 563$).

\newpage

\section{Broader utilization}

It’s not without a little luck that the sequence as described exists. If we look for $y’$ such that $y’ \cdot p\#$ rounded
to the nearest integer is prime, we wouldn’t get far: at stage 23, one prime would be left:
$\lfloor y’ \cdot 83\# +\frac12 \rfloor = q’_{23,1} = 206780313999369083332356327764879$.
There’s no prime within $q’_{23,1} \cdot 89 \pm 44$, so that sequence would be finite (or ``mortal'') --- even though the primes
are smaller compared to $y = 1.2541961...$ ($\#q’_s = \{2, 2, 2, 2, 1, 1, 2, 3, 6, 4, 1, 3, 5, 4, 2, 2,$
$2, 4, 5, 5, 2, 2, 1, 0.\}$ --- and the ratio $\frac y{y’} = 1.6198...$ is close to the golden ratio:) \\

And for those who can’t get enough sequences that grow like $y \cdot p\#$, here is a definition for ``semi-p\#Y-sequences'':
find $y$ such that $\lfloor y \cdot p\# \rfloor$ is a prime number for every $p \geq 3$. The only difference here is that $\lfloor 2y \rfloor$
is allowed to be an even number other than 2 (maybe even a semiprime, depending on the taste about restrictions, though
it’s hard to find such an example like $\lfloor y \rfloor = 210457743323$ --- with 1124 primes at $s = 331$). This leaves a lot of
sequences on the verge of extinction for quite a while during the calculations. Some of them should survive an infinite run:
\\
$0 < y < 1$: this gives a sequence which is stronger than the original p\#Y-sequence \\
($1 < y < 2$: original sequence) \\
$25 < y < 26$: stronger than the original sequence up to stage 166, weaker after stage 168 \\
$1411 < y < 1412$: a weak sequence with only 883 primes at $s = 168$, $p = 997$ \\
$3432 < y < 3433$: a very weak sequence, only 168 primes at $s = 168$, $p = 997$ \\
$13948 < y < 13949$: another weak sequence \\
$201420 < y < 201421$: another one which has somewhat good chances to survive \\
$6007103 < y < 6007104$: rather weak \\
$25510020 < y < 25510021$: just a bit stronger than $\lfloor y \rfloor = 3432$ until $s = 167$, $p = 991$ \\
\\
... Can you find others? Not like these, they’re only finite (successive records):
\\
$2 < y < 3$: extinct at stage 4 \\
$3 < y < 4$: extinct at stage 7 \\
$5 < y < 6$: extinct at stage 11 \\
$16 < y < 17$: extinct at stage 14 \\
$978 < y < 979$: extinct at stage 17 \\
$6640 < y < 6641$: extinct at stage 23 \\
$11456 < y < 11457$: extinct at stage 35 \\
$160563 < y < 160564$: extinct at stage 38 \\
$283257 < y < 283258$: extinct at stage 68 \\
$1117230 < y < 1117231$: extinct at stage 78 \\
$1594501 < y < 1594502$: extinct at stage 86 \\
$55990660 < y < 55990661$: extinct at stage 106 \\
$108286142 < y < 108286143$: extinct at stage 114 \\

The number of $y$’s with different integer part that result in infinite semi-sequences should also be infinite.

\section{Appendix A. More on the evolution}

\begin{table}[H]
\caption{Data for $s \leq 330$, using BPSW-pseudoprimality for the primes $q$ for $s \geq 110$.
$n*$ stands for the number of surviving branches, backtracked from the 810387 primes of stage 331.
This number seems to be settled for $s \leq 130$ but will be subject to change for $s \geq 131$ by the time the
calculation is taken further since many of them will get cancelled out, so this column only gives a momentary picture
(as well as an upper bound). Nevertheless, the rate of decrease in $n*$ for larger $s$ may be of interest in itself.
\label{table7}
}
\begin{tabular}{ccc|cccc|cccc|cccc}
\hline
$s$ & $p$ & $n$ & $s$ & $p$ & $n$ & $n*$ & $s$ & $p$ & $n$ & $n*$ & $s$ & $p$ & $n$ & $n*$ \\
\hline
1 & 2 & 1 & 31 & 127 & 2 & 1 & 61 & 283 & 93 & 8 & 91 & 467 & 413 & 35 \\
2 & 3 & 1 & 32 & 131 & 5 & 1 & 62 & 293 & 81 & 9 & 92 & 479 & 426 & 36 \\
3 & 5 & 1 & 33 & 137 & 6 & 1 & 63 & 307 & 81 & 9 & 93 & 487 & 446 & 37 \\
4 & 7 & 1 & 34 & 139 & 12 & 1 & 64 & 311 & 67 & 9 & 94 & 491 & 454 & 37 \\
5 & 11 & 2 & 35 & 149 & 21 & 1 & 65 & 313 & 66 & 11 & 95 & 499 & 491 & 37 \\
6 & 13 & 2 & 36 & 151 & 19 & 1 & 66 & 317 & 74 & 11 & 96 & 503 & 456 & 38 \\
7 & 17 & 3 & 37 & 157 & 15 & 1 & 67 & 331 & 91 & 12 & 97 & 509 & 490 & 38 \\
8 & 19 & 3 & 38 & 163 & 16 & 1 & 68 & 337 & 88 & 12 & 98 & 521 & 559 & 41 \\
9 & 23 & 4 & 39 & 167 & 24 & 1 & 69 & 347 & 90 & 14 & 99 & 523 & 573 & 43 \\
10 & 29 & 6 & 40 & 173 & 18 & 1 & 70 & 349 & 95 & 15 & 100 & 541 & 594 & 47 \\
11 & 31 & 4 & 41 & 179 & 18 & 1 & 71 & 353 & 102 & 16 & 101 & 547 & 652 & 47 \\
12 & 37 & 5 & 42 & 181 & 17 & 1 & 72 & 359 & 126 & 18 & 102 & 557 & 718 & 50 \\
13 & 41 & 5 & 43 & 191 & 14 & 1 & 73 & 367 & 152 & 19 & 103 & 563 & 757 & 50 \\
14 & 43 & 9 & 44 & 193 & 24 & 2 & 74 & 373 & 154 & 19 & 104 & 569 & 786 & 53 \\
15 & 47 & 11 & 45 & 197 & 24 & 3 & 75 & 379 & 166 & 19 & 105 & 571 & 835 & 55 \\
16 & 53 & 10 & 46 & 199 & 28 & 5 & 76 & 383 & 187 & 20 & 106 & 577 & 839 & 57 \\
17 & 59 & 12 & 47 & 211 & 30 & 5 & 77 & 389 & 214 & 20 & 107 & 587 & 854 & 59 \\
18 & 61 & 8 & 48 & 223 & 36 & 5 & 78 & 397 & 206 & 21 & 108 & 593 & 906 & 63 \\
19 & 67 & 6 & 49 & 227 & 49 & 5 & 79 & 401 & 201 & 21 & 109 & 599 & 916 & 65 \\
20 & 71 & 11 & 50 & 229 & 44 & 5 & 80 & 409 & 220 & 21 & 110 & 601 & 998 & 69 \\
21 & 73 & 5 & 51 & 233 & 52 & 5 & 81 & 419 & 241 & 23 & 111 & 607 & 988 & 71 \\
22 & 79 & 4 & 52 & 239 & 53 & 5 & 82 & 421 & 249 & 25 & 112 & 613 & 1016 & 74 \\
23 & 83 & 6 & 53 & 241 & 55 & 5 & 83 & 431 & 269 & 25 & 113 & 617 & 1078 & 78 \\
24 & 89 & 3 & 54 & 251 & 67 & 6 & 84 & 433 & 320 & 27 & 114 & 619 & 1165 & 81 \\
25 & 97 & 2 & 55 & 257 & 69 & 6 & 85 & 439 & 354 & 29 & 115 & 631 & 1237 & 84 \\
26 & 101 & 1 & 56 & 263 & 72 & 7 & 86 & 443 & 354 & 31 & 116 & 641 & 1295 & 86 \\
27 & 103 & 3 & 57 & 269 & 81 & 7 & 87 & 449 & 365 & 32 & 117 & 643 & 1371 & 90 \\
28 & 107 & 1 & 58 & 271 & 79 & 7 & 88 & 457 & 369 & 33 & 118 & 647 & 1399 & 93 \\
29 & 109 & 1 & 59 & 277 & 85 & 8 & 89 & 461 & 358 & 33 & 119 & 653 & 1523 & 97 \\
30 & 113 & 3 & 60 & 281 & 83 & 8 & 90 & 463 & 387 & 34 & 120 & 659 & 1618 & 101 \\
\hline
\end{tabular}
\end{table}

\newpage

Stage 139 is the last known point where the number of primes decreases. (Interestingly enough, just before that happens,
we have $n = 4p$, which is the only known --- but certainly not the only --- example where $n$ is a multiple of $p$.) \\

\begin{tabular}{cccc|cccc|cccc}
\hline
$s$ & $p$ & $n$ & $n*$ & $s$ & $p$ & $n$ & $n*$ & $s$ & $p$ & $n$ & $n*$ \\
\hline
121 & 661 & 1701 & 105 & 156 & 911 & 5971 & 387 & 191 & 1153 & 19690 & 1138 \\
122 & 673 & 1773 & 108 & 157 & 919 & 6366 & 401 & 192 & 1163 & 20364 & 1169 \\
123 & 677 & 1801 & 116 & 158 & 929 & 6582 & 411 & 193 & 1171 & 21154 & 1205 \\
124 & 683 & 1889 & 121 & 159 & 937 & 6799 & 430 & 194 & 1181 & 21954 & 1251 \\
125 & 691 & 1885 & 127 & 160 & 941 & 7006 & 452 & 195 & 1187 & 22688 & 1274 \\
126 & 701 & 1944 & 132 & 161 & 947 & 7151 & 460 & 196 & 1193 & 23528 & 1319 \\
127 & 709 & 2008 & 138 & 162 & 953 & 7599 & 471 & 197 & 1201 & 24295 & 1349 \\
128 & 719 & 2141 & 144 & 163 & 967 & 7920 & 491 & 198 & 1213 & 25085 & 1387 \\
129 & 727 & 2136 & 150 & 164 & 971 & 8276 & 508 & 199 & 1217 & 25717 & 1428 \\
130 & 733 & 2241 & 158 & 165 & 977 & 8708 & 521 & 200 & 1223 & 26680 & 1477 \\
131 & 739 & 2354 & 165 & 166 & 983 & 9086 & 532 & 201 & 1229 & 27569 & 1510 \\
132 & 743 & 2442 & 177 & 167 & 991 & 9577 & 550 & 202 & 1231 & 28242 & 1555 \\
133 & 751 & 2548 & 179 & 168 & 997 & 9736 & 566 & 203 & 1237 & 28917 & 1595 \\
134 & 757 & 2649 & 184 & 169 & 1009 & 10329 & 581 & 204 & 1249 & 29871 & 1639 \\
135 & 761 & 2791 & 193 & 170 & 1013 & 10722 & 589 & 205 & 1259 & 30707 & 1688 \\
136 & 769 & 2998 & 200 & 171 & 1019 & 11048 & 606 & 206 & 1277 & 31555 & 1733 \\
137 & 773 & 3017 & 207 & 172 & 1021 & 11387 & 634 & 207 & 1279 & 32633 & 1786 \\
138 & 787 & 3148 & 214 & 173 & 1031 & 11769 & 663 & 208 & 1283 & 33494 & 1849 \\
139 & 797 & 3132 & 224 & 174 & 1033 & 12251 & 679 & 209 & 1289 & 34556 & 1899 \\
140 & 809 & 3385 & 231 & 175 & 1039 & 12490 & 705 & 210 & 1291 & 35744 & 1957 \\
141 & 811 & 3537 & 245 & 176 & 1049 & 12905 & 727 & 211 & 1297 & 36546 & 2010 \\
142 & 821 & 3698 & 250 & 177 & 1051 & 13209 & 743 & 212 & 1301 & 37363 & 2063 \\
143 & 823 & 3914 & 261 & 178 & 1061 & 13556 & 768 & 213 & 1303 & 38018 & 2123 \\
144 & 827 & 4175 & 266 & 179 & 1063 & 13918 & 790 & 214 & 1307 & 38936 & 2188 \\
145 & 829 & 4345 & 273 & 180 & 1069 & 14393 & 815 & 215 & 1319 & 39633 & 2272 \\
146 & 839 & 4360 & 284 & 181 & 1087 & 15006 & 845 & 216 & 1321 & 40152 & 2344 \\
147 & 853 & 4537 & 291 & 182 & 1091 & 15472 & 867 & 217 & 1327 & 40822 & 2420 \\
148 & 857 & 4722 & 303 & 183 & 1093 & 15695 & 901 & 218 & 1361 & 42524 & 2500 \\
149 & 859 & 4862 & 307 & 184 & 1097 & 16075 & 930 & 219 & 1367 & 44138 & 2580 \\
150 & 863 & 4951 & 319 & 185 & 1103 & 16548 & 963 & 220 & 1373 & 45674 & 2654 \\
151 & 877 & 5099 & 331 & 186 & 1109 & 16852 & 999 & 221 & 1381 & 47308 & 2738 \\
152 & 881 & 5228 & 344 & 187 & 1117 & 17548 & 1024 & 222 & 1399 & 49203 & 2818 \\
153 & 883 & 5334 & 356 & 188 & 1123 & 17926 & 1055 & 223 & 1409 & 50999 & 2905 \\
154 & 887 & 5460 & 363 & 189 & 1129 & 18415 & 1081 & 224 & 1423 & 53578 & 2982 \\
155 & 907 & 5715 & 377 & 190 & 1151 & 19145 & 1115 & 225 & 1427 & 56009 & 3071 \\
\hline
\end{tabular}

\newpage

\begin{tabular}{cccc|cccc|cccc}
\hline
$s$ & $p$ & $n$ & $n*$ & $s$ & $p$ & $n$ & $n*$ & $s$ & $p$ & $n$ & $n*$ \\
\hline
226 & 1429 & 58245 & 3143 & 261 & 1663 & 141747 & 8077 & 296 & 1949 & 323493 & 26835 \\
227 & 1433 & 60614 & 3264 & 262 & 1667 & 144129 & 8290 & 297 & 1951 & 332369 & 28016 \\
228 & 1439 & 63012 & 3359 & 263 & 1669 & 146278 & 8514 & 298 & 1973 & 344077 & 29288 \\
229 & 1447 & 65622 & 3434 & 264 & 1693 & 150011 & 8752 & 299 & 1979 & 355840 & 30683 \\
230 & 1451 & 67609 & 3542 & 265 & 1697 & 153719 & 9016 & 300 & 1987 & 367371 & 32050 \\
231 & 1453 & 69411 & 3650 & 266 & 1699 & 157438 & 9268 & 301 & 1993 & 379303 & 33504 \\
232 & 1459 & 71677 & 3745 & 267 & 1709 & 160935 & 9582 & 302 & 1997 & 391133 & 35062 \\
233 & 1471 & 73978 & 3854 & 268 & 1721 & 164713 & 9877 & 303 & 1999 & 402306 & 36665 \\
234 & 1481 & 76816 & 3964 & 269 & 1723 & 168209 & 10187 & 304 & 2003 & 413306 & 38433 \\
235 & 1483 & 78753 & 4049 & 270 & 1733 & 171914 & 10497 & 305 & 2011 & 424348 & 40351 \\
236 & 1487 & 81051 & 4156 & 271 & 1741 & 175208 & 10838 & 306 & 2017 & 435084 & 42411 \\
237 & 1489 & 82647 & 4260 & 272 & 1747 & 178924 & 11177 & 307 & 2027 & 446383 & 44703 \\
238 & 1493 & 84106 & 4386 & 273 & 1753 & 182442 & 11549 & 308 & 2029 & 457067 & 47035 \\
239 & 1499 & 85756 & 4511 & 274 & 1759 & 186286 & 11951 & 309 & 2039 & 468311 & 49593 \\
240 & 1511 & 87469 & 4631 & 275 & 1777 & 190466 & 12354 & 310 & 2053 & 482068 & 52470 \\
241 & 1523 & 90333 & 4765 & 276 & 1783 & 196294 & 12773 & 311 & 2063 & 494789 & 55576 \\
242 & 1531 & 93284 & 4879 & 277 & 1787 & 200525 & 13184 & 312 & 2069 & 509095 & 58953 \\
243 & 1543 & 96279 & 5012 & 278 & 1789 & 204944 & 13625 & 313 & 2081 & 524247 & 62650 \\
244 & 1549 & 99150 & 5141 & 279 & 1801 & 209420 & 14097 & 314 & 2083 & 538695 & 66722 \\
245 & 1553 & 101723 & 5281 & 280 & 1811 & 214672 & 14574 & 315 & 2087 & 551951 & 71192 \\
246 & 1559 & 104560 & 5402 & 281 & 1823 & 218990 & 15100 & 316 & 2089 & 564045 & 76342 \\
247 & 1567 & 107041 & 5542 & 282 & 1831 & 224272 & 15623 & 317 & 2099 & 577252 & 81969 \\
248 & 1571 & 109795 & 5691 & 283 & 1847 & 230486 & 16200 & 318 & 2111 & 592642 & 88171 \\
249 & 1579 & 111912 & 5827 & 284 & 1861 & 238339 & 16792 & 319 & 2113 & 605293 & 95404 \\
250 & 1583 & 114508 & 5989 & 285 & 1867 & 246127 & 17405 & 320 & 2129 & 620280 & 103655 \\
251 & 1597 & 117842 & 6160 & 286 & 1871 & 253236 & 18037 & 321 & 2131 & 635284 & 113441 \\
252 & 1601 & 121087 & 6327 & 287 & 1873 & 259946 & 18726 & 322 & 2137 & 649598 & 124938 \\
253 & 1607 & 124407 & 6493 & 288 & 1877 & 266707 & 19481 & 323 & 2141 & 664112 & 138630 \\
254 & 1609 & 126382 & 6682 & 289 & 1879 & 272114 & 20242 & 324 & 2143 & 676043 & 155455 \\
255 & 1613 & 128949 & 6854 & 290 & 1889 & 278512 & 21025 & 325 & 2153 & 689579 & 176551 \\
256 & 1619 & 130732 & 7021 & 291 & 1901 & 285821 & 21851 & 326 & 2161 & 703065 & 203256 \\
257 & 1621 & 132088 & 7226 & 292 & 1907 & 292121 & 22714 & 327 & 2179 & 720414 & 240371 \\
258 & 1627 & 133678 & 7440 & 293 & 1913 & 298404 & 23666 & 328 & 2203 & 742364 & 292120 \\
259 & 1637 & 135308 & 7642 & 294 & 1931 & 305779 & 24693 & 329 & 2207 & 763990 & 370783 \\
260 & 1657 & 137970 & 7862 & 295 & 1933 & 314060 & 25728 & 330 & 2213 & 787038 & 506979 \\
\hline
\end{tabular}

\vspace{8mm}

Data for $s=331$ as output from the program on the following page is given as an ancillary file to this version of the paper.

\newpage

The calculation up to the point given above takes about eight months with Pari/GP \cite{Pari} on a
single CPU core of a 2020 state-of-the-art PC with the following program, which in the
given form starts at $s$ = 29 and keeps track of the numbers as {\tt a} = $q_{s,1}$ and a memory
friendly vector {\tt d} comprising only the consecutive differences between the numbers of a
given stage where {\tt d[1]} = 0 and, for $i > 1$, {\tt d[i]} = $q_{s,i}-q_{s,i-1}$:
\footnote{~{\tt a} and {\tt d} can also be read in from a previously calculated output file} \\

\begin{verbatim}
{
a=350842542483891235293716663559065020274899073; d=[0];

s=0; b=a;
while(b>1, s++; b\=prime(s));
i=#d; gettime();
while(s<330,
   s++; p=prime(s); o=a*p; c=d;
   e=99+floor(i*(1+2/sqrt(p))); d=vector(e);
   m=i; i=0;
   for(j=1, m,
      print1("interval "j"/"m,Strchr(13));
      o+=c[j]*p; v=vector(p);
      forprime(b=3, p-2, r=b-lift(Mod(o,b)); forstep(l=r, p, b, v[l]=1));
      forprime(b=p+2, p^2, r=b-lift(Mod(o,b)); if(r<p, v[r]=1));
      forstep(r=2, p-1, 2,
         if(!v[r],
            q=o+r;
            if(ispseudoprime(q),
               i++; if(i>1, d[i]=q-z, a=q); z=q;
               print1("interval "j"/"m" - "i" primes found",Strchr(13))
            )
         )
      )
   );
   g=floor(gettime()/1000); x="[CPU time: ";
   f=floor(g/3600); if(f,x=Str(x,f"h "));
   f=floor(g/60); if(f,x=Str(x,f%60"m "));
   x=Str(x,g%60"s] "); t=Str("Stage "s); print(t" (p="p"): "i" primes "x);
   t=Str("p#Y "t".txt"); write(t,"a="a"; d="vecextract(d,Str("1.."i)))
)
}

\end{verbatim}

Be sure to set {\tt allocatemem(2*10$^\wedge$8)} or higher (and watch for the second ``while'' condition)
when aiming for $s > 330$.

\newpage

\section{Appendix B. ``Mille''stone}

The first titanic p\#Y-prime ($s = 350$, $p = 2357$, primality proved over all stages): $q_{350,1} =$
297588802944669004567432988935086641892218253950808346654606799666416703403235386
110626587476071000262691671679224513211575176049294913346932056675348621165796730
263895561614248056123980950034528686832211500014435248604833958209819278373161550
816949962259892736268616001731340442751214329425997049930689376917146879956304444
951698385088539246081531976376843432978119495609784805439847386853846175609773406
294780298954522227295543214788475948953215376258100904758911480647631719489690751
598932601226480640025811451192259365414314027908668601256915767355990441760663317
227840259709936806311906508957316674862617637555729680550552499945756091942766728
462111799759803946513079633400975724089109421146044040340589262506167584698451951
787089619451193713787721236181896575926339065362406372104874478640876833784201015
302057350287221183373441732847205339141787284034897239162684632744754226872565733
662978724742861772240754047058150049733225603812946433723077381838542363028349647
7698386126623480539308646019 \\
It’s a nice fact that 2357\# of all primorials has exactly 1000 digits. \\

The ostensibly first titanic p\#Y-prime that survives in the long run, the $16^{th}$ smallest of
them all at $p = 2357$, spanning its branch beyond 13751\# ($s > 1626$), is $q_{350,1} +$
648890029246175571592953681722255270285895737015173251749453742478642415683875459
622820131990229693513724200589337798081947713889805053628824480433575300758724935
31089513889637480426037406091083077809690, \\
which leads to $y_{min} \geq$ \\
1.25419610157801193627767955491421342377986921804262219583272255460886469942875144
751323169673647331200713029313835829410519055068071464454595777347989721472473549
850870383774045538648544910432104569800625071405146132797606019221734399136669410
231685026270656941987419822020206973004280705089121259525801855611303610702614940
065145204852907234077802800854431673412825358322402127785595190344370611325045336
730133684806271064030484811296434802163878424799778915265485948593800601469382480
980040875271653561997759928697347239814789283725106267672399076512871960330242356
126907949160683090804262705523408008158909911209586208398778315843260104151172880
387141761635412295558476482556483533650156157283843179063498926572476414145174204
072485693292701787200666665412458613015446492323282769411572687097469909612014411
743405605429162630100500385105709425493894716245564452834943026368056855139809017
404605931009796839371601349304261945450753825556512381095336296655187077486338594
1160720014897114954349144158... \\

The first 129 primes of stage 428 ($p = 2969$) are known to be ``mortal'': $q_{428,130}$ is the first
that may survive an infinite run, which is a currently known maximum. It takes on average
about $\sqrt{p}/2$ primes to find a surviving prime. An error term for this approximation,
assumed to be of order of at most $O(\sqrt{p}/\log p)$, would probably be desireable.

\section{Appendix C. Sequence tuples (multiple primes in one interval)}
$ $\\
First twin: ($1^{st}$ prime of stage 4)*11+ \{4, 10\} \\
First triplet: ($1^{st}$ prime of $s$ = 7)*19+ \{10, 12, 16\} \\
First quadruplet: ($2^{nd}$ prime of $s$ = 9)*29+ \{6, 14, 24, 26\} \\
First quintuplet: ($18^{th}$ prime of $s$ = 46)*211+ \{42, 110, 140, 144, 194\} \\
First sextuplet: ($69^{th}$ prime of $s$ = 69)*349+ \{18, 40, 234, 262, 292, 298\} \\
First septuplet: ($238^{th}$ prime of $s$ = 84)*439+ \{34, 228, 252, 282, 364, 378, 382\} \\
First octuplet: ($5687^{th}$ prime of $s$ = 159)*941+ \{170, 234, 294, 462, 696, 740, 752, 812\} \\
First nonuplet: ($270^{th}$ prime of $s$ = 181)*1091+ \{18, 46, 60, 166, 180, 312, 850, 1062, 1080\} \\
First decuplet: ($193057^{th}$ prime of $s$ = 289)*1889 + \{220, 238, 378, 624, 934, 1048, 1414, \\ 
\phantom{.} \hfill 1612, 1678, 1750\} \\
\\
Higher-order descendants (first instances of number of ``grandchildren''): \\
3 descendants of order 2: ($1^{st}$ prime of stage 6)*17*19+ \{162, 164, 168\} \\
4 descendants: ($1^{st}$ prime of $s$ = 7)*19*23+ \{246, 284, 384, 386\} \\
5 descendants: ($5^{th}$ prime of $s$ = 21)*79*83+ \{3504, 3520, 5190, 5200, 5224\} \\
6 descendants: ($4^{th}$ prime of $s$ = 12)*43*47+ \{102, 108, 582, 598, 1428, 1450\} \\
7 descendants: ($7^{th}$ prime of $s$ = 34)*149*151+ \{754, ..., 12484\} \\
8 descendants: ($30^{th}$ prime of $s$ = 55)*263*269+ \{7686, ..., 40560\} \\
9 descendants: ($10^{th}$ prime of $s$ = 59)*281*283+ \{11388, ..., 59016\} \\
10 descendants: ($16^{th}$ prime of $s$ = 47)*223*227+ \{13650, ..., 38340\} \\
11 descendants: ($42^{nd}$ prime of $s$ = 67)*337*347+ \{24432, ..., 116892\} \\
12 descendants: ($27^{th}$ prime of $s$ = 55)*263*269+ \{45742, ..., 68590\} \\
13 descendants: ($3053^{rd}$ prime of $s$ = 138)*797*809+ \{195918, ..., 603780\} \\
14 descendants: ($623^{rd}$ prime of $s$ = 101)*557*563+ \{94812, ..., 261562\} \\
15 descendants: ($32^{nd}$ prime of $s$ = 72)*367*373+ \{2424, ..., 134612\} \\
16 descendants: ($45231^{st}$ prime of $s$ = 224)*1427*1429+ \{74772, ..., 1944498\} \\
17 descendants: ($2030^{th}$ prime of $s$ = 176)*1051*1061+ \{74496, ..., 1083072\} \\
18 descendants: ($376869^{th}$ prime of $s$ = 321)*2137*2141+ \{1058496, ..., 3727156\} \\

The like-minded reader with enough time on their hands might be encouraged to find successive
maxima of descendants of higher order. To start, here’s 4 descendants of order 3: \\
($1^{st}$ prime of stage 6)*17*19*23+ \{3742, 3780, 3880, 3882\}

\newpage

\section{Appendix D. Relative strength of split intervals}

\begin{figure}[H]
  \centering
  \includegraphics[bb=-40 0 840 560,width=8in,height=5.5in,keepaspectratio]{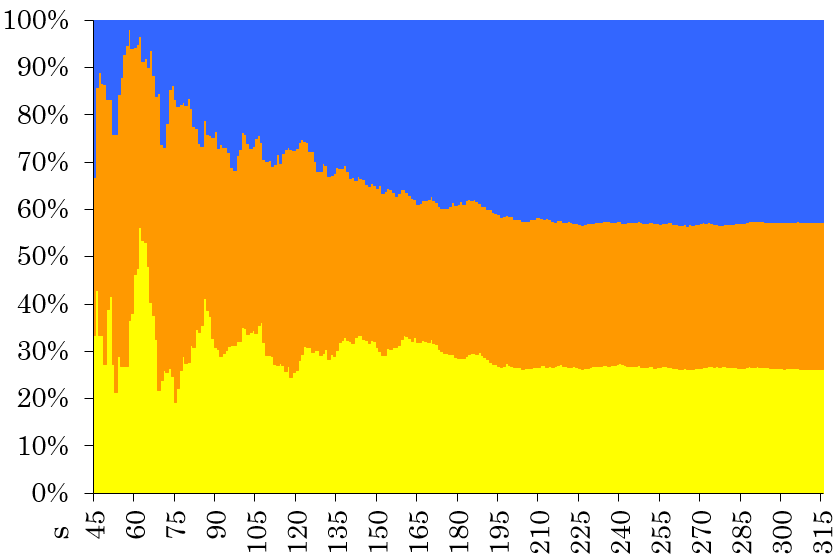}
  \caption{A nice and colorful edit showing the strength of the first three split intervals, starting
at stage 45. The 592642 probable primes at stage 318 are here divided into three groups
(154145 + 184653 + 253844) with the respective relative strength 26.01\%, 31.16\%, and
42.83\%. These values may (hopefully) stabilize while the number of primes in the sequence
continues to grow. Giving the exact percentages as $s$ goes to infinity will be rather difficult.}
\label{fig6}
\end{figure}

\newpage

\section{Appendix E. Plouffe, revisited}

As a final tidbit of numerical info, inspired by Plouffe's work mentioned in the introduction, we give a
constant $A$ (to 1029 places after the decimal point) such that $\lfloor A^{1.001^n} \rfloor$ should
give an infinite sequence of primes, for integer $n \geq 0$: \\

$A = 10^{3875} -$ 9840440.013269045402699925892590230034000537023543058287831784100957
276330984498140336774521413625639037370179638123065989644921831945828284474797146
214314469465943653321876212421217357462243438103894210660374222702124802145689793
662592212239695486359151401226648141442872210231891779023391189805301737658502412
189281932141114052020391866640443434821779837126175086985360497145764015143114697
436609454475477982979703079133332121473465799127802190473213691503591043455180286
502985352685402044844822486579962417722857237662143660192493843854713085865446265
026287773151577320009216914265110643337790882410603483306933710003501741787035605
297213202524866787149107203669351182374236682986886760407222192912257430323369703
769706007229005048714432375371543157174405792236627110947314049230460840453363435
080982900105891219946283712237437066606879315185506382375652998319912097105539735
706790752939943728420192710162069587039287823862485223550170317760285385396844814
481707064610195932400139047255566802623529101069361896409968571400144014982338... \\

We would be happy to learn about a smaller known value for $A$.

\vspace{20mm}

\section{Acknowledgements}

Many thanks to Charles R. Greathouse IV for his support when we first presented this sequence
at {\tt Mersenneforum.org} back in 2008. We are also deeply grateful to Alexei Kourbatov for helpful
suggestions on the paper.

\vspace{15mm}
{\small\bigskip\noindent
{\bf Keywords:} 
prime number sequence, prime generating formula, primorial function, Chebyshev's bias.
}

\vspace{15mm}
{\small\bigskip\noindent

{\bf About the author:}
Martin Raab works in computational number theory as a hobby. He performed various computations regarding
prime numbers, some of which are found in the On-Line Encyclopedia of Integer Sequences (\url {https://oeis.org/}),
at the Mersenneforum (\url {https://www.mersenneforum.org/index.php}), and in the tables of first occurrence prime
gaps (currently hosted at \url {https://primegap-list-project.github.io/}). Martin lives in 63840 Hausen, Germany.
He can be reached at kilroy14159265@gmail.com.
}

\end{document}